\documentclass[12pt, reqno]{amsart}
\usepackage{amsmath, amsthm, amscd, amsfonts, amssymb, graphicx, color}
\usepackage[bookmarksnumbered, colorlinks, plainpages]{hyperref}
\hypersetup{colorlinks=true,linkcolor=red, anchorcolor=green, citecolor=cyan, urlcolor=red, filecolor=magenta, pdftoolbar=true}

\newtheorem{theorem}{Theorem}[section]

\theoremstyle{definition}

\theoremstyle{remark}

\numberwithin{equation}{section}

\begin{document}

\setcounter{page}{1}

\title[Estimates for Mittag-Leffler functions with smooth phase depending on two variables]{Estimates for Mittag-Leffler functions with smooth phase depending on two variables}

\author[A.R.Safarov]{Akbar R.Safarov}

\address{\textcolor[rgb]{0.00,0.00,0.84}{Akbar R.Safarov \newline Uzbekistan Academy of Sciences V.I.Romanovskiy Institute of Mathematics \newline Olmazor district, University 46, Tashkent, Uzbekistan \newline Samarkand State University \newline Department Mathematics, \\
15 University Boulevard
 \newline  Samarkand, 140104, Uzbekistan }}
\email{\textcolor[rgb]{0.00,0.00,0.84}{safarov-akbar@mail.ru}}





\begin{abstract}
In this paper  we  consider the problem on estimates for Mittag-Leffler functions with the smooth phase functions of two variables having singularities of type $D_{\infty} $, $D_{4}^{\pm}$ and $A_{r}$. The generalisation is that we replace the exponential function with the Mittag-Leffler-type function, to study oscillatory type integrals. We extend results of paper \cite{Ruzhansky} and \cite{Ruzhansky2021} to two-dimensional integrals with phase having some simple singularities.
\end{abstract}
\maketitle
\tableofcontents
 \section{Introduction}

The Mittag-Leffler function $E_{\alpha}(z)$ is named after the great Swedish mathematician G\"{o}sta Magnus Mittag-Leffler (1846-1927) who defined it by a power series
\begin{equation}\label{Formul3}
E_{\alpha}(z)=\sum_{k=0}^{\infty}\frac{z^{k}}{\Gamma(\alpha k+1)},\,\,\ \alpha\in\mathbb{C},\,\, Re(\alpha)>0,
\end{equation}
and studied its properties in 1902-1905 in five subsequent notes \cite{ML1}-\cite{ML4} in connection with his summation method for divergent series.

 A classic generalization of the Mittag-Leffler function, namely the two-parametric Mittag-Leffler function
\begin{equation}\label{Formul4}
E_{\alpha,\beta}(z)=\sum_{k=0}^{\infty}\frac{z^{k}}{\Gamma(\alpha k+\beta)},\,\,\ \alpha,\beta\in\mathbb{C},\,\, Re(\alpha)>0,
\end{equation}
which was deeply investigated independently by Humbert and Agarval in 1953 \cite{Hum53}-\cite{HumAga53} and by Dzherbashyan in 1954
\cite{Dzh54a} - \cite{Dzh54c}, see also  \cite{Rudolf} and the references therein.

 In this paper we also consider a special case   the generalized Mittag-Leffler function defined as in \eqref{Formul4} by

$$E_{\alpha,\beta}(x)=\sum_{k=0}^{\infty}\frac{x^{k}}{\Gamma(\alpha k+\beta)}, \alpha>0, \beta\in\mathbb{R}.$$
Obviously,
\begin{equation}\label{Form1}
E_{1,1}(x)=e^{x}.
\end{equation}
We consider the following integral with phase $f$ and amplitude $\psi$
\begin{equation}\label{int1}
I_{\alpha,\beta}=\int_{a}^{b} E_{\alpha,\beta}(i\lambda f(x))\psi(x)dx
\end{equation}
where $0<\alpha\leq1,$ $\beta>0$ and $\lambda>0$.

If $\alpha=\beta=1$ in the integral \eqref{int1}. The  integral \,$I_{1,1}$ \,is called an oscillatory integral. In harmonic analysis estimates, the most important estimates for oscillatory integrals van der Corput lemma \cite{VanDer}.  Estimates for oscillatory integrals with polynomial phase can be bound, for instance, in papers  \cite{AKC}-\cite{Safarov1}. In the current paper we replace the exponential function with the Mittag-Leffler-type function and study oscillatory type integrals \eqref{int1}. In the papers \cite{Ruzhansky} and \cite{Ruzhansky2021}  analogues of the van der Corpute lemmas involving Mittag-Leffler functions for one dimensional integrals have been considered. We extend results of \cite{Ruzhansky} and \cite{Ruzhansky2021} for two-dimensional integrals with phase having some simple singularities.

  The main result of the paper is the following.

\begin{theorem}\label{Theor1}
Let $-\infty<a<b<\infty$. Assume that the phase function is a homogenous  polynomial of third degree with two variables and let $\psi\in L^{p}[a,b]^{2}$, $1<p\leq\infty$.
Then for any $\alpha\in(0,1)$, $\beta,\lambda\in(0,+\infty)$
\begin{align}\label{Formul5}
& \left|\int_{[a,b]^2}E_{\alpha,\beta}(i\lambda x_{1}^{2}x_{2})\psi(x)dx\right|\leq\frac{C\|\psi\|_{L^{p}}}{\lambda^{\frac{1}{2}-\frac{1}{2p}}}, \\
& \left|\int_{[a,b]^2}E_{\alpha,\beta}(i\lambda (x_{1}^{2}x_{2}\pm x_{2}^{3}))\psi(x)dx\right|\leq\frac{C\|\psi\|_{L^{p}}}{\lambda^{\frac{2}{3}-\frac{1}{3p}}}, \\
&\left|\int_{[a,b]^2}E_{\alpha,\beta}(i\lambda x_{1}^{3})\psi(x)dx\right|\leq\frac{C\|\psi\|_{L^{p}}}{\lambda^{\frac{1}{3}-\frac{1}{3p}}},
\end{align}
where constant $C$ depends only on $p.$
\end{theorem}

\section{Some auxiliary statements}
First we give  ancillary statements. Let us consider homogenous  polynomial third degree with two variables.

\textbf{Proposition 1 }\label{Prop1}(\cite{Arn}).
A homogenous  polynomial of third degree with two variables may be reduced by a R-linear transformation to one of the forms:
1) $x_{1}^{2}x_{2}$, 2) $x_{1}^{2}x_{2}\pm x_{2}^{3}$, 3) $x_{1}^{3}$, 4) 0.

\textbf{Definition 1.}
Given $\mu\in(1,\infty]$, a critical point, equivalent to the critical point of the function
 $x_{1}^{2}x_{2}\pm x_{2}^{\mu-1}$ is said to be a critical point of type $D_{\mu}^{\pm}$,  where $x_{2}^{\mu-1}\equiv0$ for $\mu=\infty$.

\textbf{Definition 2.}
A critical point, equivalent to the critical point of the function
 $x_{1}^{r+1},$ $r\geq1$ is said to be a critical point of type $A_{r}$.

\textbf{Proposition 2 }(\cite{Pod}).
If $0<\alpha<2,\beta$ is an arbitrary real number and $\mu$ is such that $\pi\alpha/2<\mu<\min\{\pi,\pi\alpha\},$ then there is $C>0$ such that
\begin{equation}\label{Formul1}
  |E_{\alpha,\beta}(z)|\leq\frac{C}{1+|z|}, z\in\mathbb{C},\mu\leq|\arg(z)|\leq\pi.
\end{equation}

\textbf{Proposition 3.}\cite{Ruzhansky}
Let $\alpha,\beta>0$ and $f:[a,b]\rightarrow\mathbb{C}$. Then for all $\lambda\in\mathbb{C}$
\begin{equation}\label{Formul2}
E_{\alpha,\beta}(i\lambda f(x))=E_{2\alpha,\beta}(-\lambda^{2}f^{2}(x))+i\lambda f(x)E_{2\alpha,\beta+\alpha}(-\lambda^{2}f^{2}(x)).
\end{equation}

\section{Proof of the main results}

\textbf{Proof of theorem \ref{Theor1}.}
As for small $\lambda$ the integral \eqref{int1} is clearly bounded, we prove Theorem \ref{Theor1} only for $\lambda\geq1.$ Without loss of generality, we can consider the integral on $[0,1]^{2}$, otherwise we reduce to this case using a linear transformation.
Since we are given a homogenous  polynomial of third degree with two variables, Proposition 1 we can represent it one of the following: 1) $x_{1}^{2}x_{2}$; 2) $x_{1}^{2}x_{2}\pm x_{2}^{3}$; 3) $x_{1}^{3}$ 4) 0. If phase function is $f(x)\equiv0$ it is clear that integral will be identically zero. So we will consider the other three cases separately.

Using inequalities \eqref{Formul1} and \eqref{Formul2} we obtain:
$$|E_{\alpha,\beta}(i\lambda f(x))|\leq|E_{2\alpha,\beta}(-\lambda^{2}f^{2}(x))|+\lambda |f(x)||E_{2\alpha,\beta+\alpha}(-\lambda^{2}f^{2}(x))|$$
\begin{equation}\label{int10}\leq\frac{C}{1+\lambda^{2}f^{2}(x)}+\frac{C\lambda|f(x)|}{1+\lambda^{2}f^{2}(x)}\leq\frac{C(1+\lambda|f(x)|)}{1+\lambda^{2}f^{2}(x)}\leq \frac{C}{1+\lambda|f(x)|}.
\end{equation}

\textbf{Case I.} First  we assume that the phase function has critical point of type $D_{\infty}$ so that $f(x)=x_{1}^{2}x_{2}.$

We consider integral \eqref{int1} of the form:

\begin{equation}\label{int16}
I_{\alpha,\beta}=\int _{[0,1]^{2} }E_{\alpha,\beta}(i\lambda x_{1}^{2}x_{2})\psi(x)dx.
\end{equation}

We use inequality \eqref{int10} in the integral \eqref{int16} and we obtain:
\begin{align}
|I_{\alpha,\beta}|=& \left|\int _{[0,1]^{2} }E_{\alpha,\beta}(i\lambda x_{1}^{2}x_{2})\psi(x)dx\right|\leq\int _{[0,1]^{2} }\left|E_{\alpha,\beta}(i\lambda x_{1}^{2}x_{2})\right|\left|\psi(x)\right|dx\nonumber
 \\
\leq & C  \int_{0}^{1}dx_{1}\int_{0}^{1}\frac{|\psi(x)|dx_{2}}{1+\lambda x_{1}^{2}x_{2}}. \label{int_estimaqqq}
\end{align}

Let  $q$ be such that $\frac{1}{p}+\frac{1}{q}=1.$   Assume first that $p\neq\infty,$ so that $q>1.$ Then using the  H\"older inequality in the inner integral  we get
\begin{align*}
 J_{in1}:=  & \int_{0}^{1}\frac{|\psi(x)|dx_{2}}{1+\lambda x_{1}^{2}x_{2}}
\leq
\left(\int_{0}^{1}|\psi(x)|^{p}dx_{2}\right)^{\frac{1}{p}}\left(\int_{0}^{1}\frac{dx_{2}}{|1+\lambda x_{1}^{2}x_{2}|^{q}}\right)^{\frac{1}{q}}\\
=&
\left(\int_{0}^{1}|\psi(x)|^{p}dx_{2}\right)^{\frac{1}{p}}\left(\frac{1 - (1 + \lambda x_1^2)^{1-q}}{(q-1)\lambda x_1^2} \right)^{\frac{1}{q}}.
\end{align*}

Thus,
\begin{gather*}
|I_{\alpha,\beta}| \le C\int_0^1 \left(\int_{0}^{1}|\psi(x)|^{p}dx_{2}\right)^{\frac{1}{p}}\left(\frac{1 - (1 + \lambda x_1^2)^{1-q}}{(q-1)\lambda x_1^2} \right)^{\frac{1}{q}} dx_1.
\end{gather*}
Then using again the  H\"older  inequality in  this integral we obtain
\begin{gather*}
|I_{\alpha,\beta}|
%
\leq  C \left(\int_{0}^{1}\int_{0}^{1}|\psi(x)|^{p}dx_{2}dx_{1}\right)^{\frac{1}{p}}\left(\int_{0}^{1}\frac{1 - (1+\lambda x_1^2)^{1-q}}{(q-1) \lambda x_1^2}\,dx_1\right)^{\frac{1}{q}}\\
%
%
 \leq C\|\psi\|_{L^{p}}\left(\int_{0}^{1}\frac{1 - (1+\lambda x_1^2)^{1-q}}{(q-1) \lambda x_1^2}\,dx_1\right)^{\frac{1}{q}}.
\end{gather*}

Let
$$
K:=\int_{0}^{1}\frac{1 - (1+\lambda x_1^2)^{1-q}}{(q-1) \lambda x_1^2}\,dx_1.
$$

Since $(1+\lambda x_1^2)^{1-q}=1 + O(\lambda x_1^2)$ near $x_1=0$ and $q>1,$ the integral $K$ is convergent. To estimate $K$, first we use the change of variables $t=\sqrt\lambda x_1$ to get
\begin{align*}
K= & \frac{1}{(q-1)\sqrt\lambda} \int_0^{\sqrt\lambda} \frac{1 - (1+t^2)^{1-q}}{t^2}dt=
\frac{1}{(q-1)\sqrt\lambda} \int_0^{\sqrt\lambda} \frac{(1+t^2)^{q-1}-1}{t^2(1+t^2)^{q-1}}dt\\
= & \frac{1}{(q-1)\sqrt\lambda} \int_0^{1} \frac{(1+t^2)^{q-1}-1}{t^2(1+t^2)^{q-1}}dt +
\frac{1}{(q-1)\sqrt\lambda} \int_1^{\sqrt\lambda} \frac{(1+t^2)^{q-1}-1}{t^2(1+t^2)^{q-1}}dt =:K_1 + K_2.
\end{align*}
Since $q-1\le [q],$ where $[q]\ge1$ is the integer part of $q>1,$ by the Newton's binomial formula
$$
(1+t^2)^{q-1} \le (1+t^2)^{[q]} =1 + [q]t^2 + \frac{[q]([q]-1)}{2}t^4 + \ldots + t^{2[q]},
$$
and hence
$$
K_1=\frac{1}{(q-1)\sqrt\lambda} \int_0^{1} \frac{(1+t^2)^{q-1}-1}{t^2(1+t^2)^{q-1}}dt \le \frac{C_q}{\sqrt\lambda},
$$
where
$$
C_q:= \frac{1}{q-1}\int_0^1\frac{[q] + \frac{[q]([q]-1)}{2}t^2 + \ldots + t^{2[q]-2}}{(1+t^2)^{q-1}}dt.
$$
Moreover, since
$
\frac{(1+t^2)^{q-1}-1}{t^2(1+t^2)^{q-1}} <\frac{1}{t^2},
$
\begin{align*}
K_2= &\frac{1}{(q-1)\sqrt\lambda} \int_1^{\sqrt\lambda} \frac{(1+t^2)^{q-1}-1}{t^2(1+t^2)^{q-1}}dt <
\frac{1}{(q-1)\sqrt\lambda} \int_1^{\sqrt\lambda} \frac{1}{t^2}dt\\
= &
\frac{1}{(q-1)\sqrt\lambda}\Big(1 - \frac{1}{\sqrt{\lambda}}\Big)
<\frac{1}{(q-1)\sqrt{\lambda}}.
\end{align*}
Hence,
$$
K\le \frac{C_q'}{\sqrt{\lambda}},\qquad C_q': =C_q +\frac{1}{q-1},
$$
and
$$
|I_{\alpha,\beta}| \le \frac{C_q''\|\psi\|_{L^p}}{\lambda^{\frac{1}{2q}}},
$$
where $C_q''$ is some coefficient depending only on $q$, and hence only on $p$.

Now we consider the case $q=1.$  Notice that the coefficient $C_q''\to+\infty$ as $q\to1$ and therefore we cannot directly conclude the requires estimate from the one for $q>1$. As $q=1$, we have $p=\infty$ and $\psi\in\mathbb{L}^{\infty}.$  In view of \eqref{int_estimaqqq}, first  we estimate inner integral as
$$|J_{in1}|=\int_{0}^{1}\frac{|\psi(x)|dx_{2}}{1+\lambda x_{1}^{2}x_{2}}\leq\sup\limits_{x_{2}\in[0,1]}|\psi(x)|\int_{0}^{1}\frac{dx_{2}}{1+\lambda x_{1}^{2}x_{2}}\leq$$$$\leq\left.\frac{\sup\limits_{x_{2}\in[0,1]}|\psi(x)|}{\lambda x_{1}^{2}}\ln(1+\lambda x_{1}^{2}x_{2})\right|_{0}^{1}
=\frac{\sup\limits_{x_{2}\in[0,1]}|\psi(x)|\ln(1+\lambda x_{1}^{2})}{\lambda x_{1}^{2}}.$$

Thus
$$|I_{\alpha,\beta}|\leq\int_{0}^{1}\frac{\sup\limits_{x_{2}\in[0,1]}|\psi(x)|\ln(1+\lambda x_{1}^{2})}{\lambda x_{1}^{2}}dx_{1}\leq
C\|\psi\|_{L^{\infty}}\int_{0}^{1}\frac{\ln(1+\lambda x_{1}^{2})}{\lambda x_{1}^{2}}dx_{1}.$$
We use change variables as $\lambda x_{1}^{2}=y$ in the last integral and get
\begin{equation}\label{F55}
|I_{\alpha,\beta}|\leq\frac{C\|\psi\|_{L^{\infty}}}{\lambda^{\frac{1}{2}}}\int_{0}^{\lambda}\frac{\ln(1+y)}{y^{\frac{3}{2}}}dy\leq
\frac{C\|\psi\|_{L^{\infty}}}{\lambda^{\frac{1}{2}}}\int_{0}^{\infty}\frac{\ln(1+y)}{y^{\frac{3}{2}}}dy.
\end{equation}
Note that the last integral converges. Instead, using integration by parts  we obtain
$$\int_{0}^{\infty}\frac{\ln(1+y)}{y^{\frac{3}{2}}}dy=-\left.\lim_{N_{1}\rightarrow0,\\N_{2}\rightarrow\infty}\frac{2\ln(1+y)}{y^{\frac{1}{2}}}\right|_{N_1}^{N_2}+
\int_{0}^{\infty}\frac{2dy}{(1+y)y^{\frac{1}{2}}}=$$$$=\int_{0}^{\infty}\frac{4dy^{\frac{1}{2}}}{1+y}=\left.4\arctan y\right|_{0}^{\infty}=2\pi.$$

Thus from \eqref{F55} we get
$$|I_{\alpha,\beta}|\leq\frac{C\|\psi\|_{L^{\infty}}}{\lambda^{\frac{1}{2}}}.$$

\textbf{Case II.} Assume that the phase function has critical point of type $D_{4}^{\pm}$ so that
$f(x)=x_{1}^{2}x_{2}\pm x_{2}^{3}$. We estimate integral \eqref{int1} when the phase function has critical point of type $D_{4}^{+}$ and the case  $D_{4}^{-}$ can be done similarly.

We consider the integral

\begin{equation}\label{int17}
I_{\alpha,\beta}=\int _{[0,1]^{2} }E_{\alpha,\beta}(i\lambda(x_{1}^{2}x_{2}+x_{2}^{3}))\psi(x)dx.
\end{equation}

Using inequality \eqref{int10} for the integral \eqref{int17} we get

$$|I_{\alpha,\beta}|=\left|\int _{[0,1]^{2} }E_{\alpha,\beta}(i\lambda(x_{1}^{2}x_{2}+x_{2}^{3}))\psi(x)dx\right|\leq\int _{[0,1]^{2} }\left|E_{\alpha,\beta}(i\lambda (x_{1}^{2}x_{2}+x_{2}^{3}))\right|\left|\psi(x)\right|dx$$$$\leq\int_{0}^{1}dx_{1}\int_{0}^{1}\frac{|\psi(x)|dx_{1}}{1+\lambda(x_{1}^{2}x_{2}+x_{2}^{3})}=
\int_{0}^{1}dx_{1}\int_{0}^{1}\frac{|\psi(x)|dx_{1}}{1+\lambda x_{2}^{3}+\lambda x_{1}^{2}x_{2}}.$$
We use H\"older inequality for the last inner integral we obtain

$$|J_{in2}|:=\int_{0}^{1}\frac{|\psi(x)|dx_{1}}{|1+\lambda x_{2}^{3}+\lambda x_{2}x_{1}^{2}|}\leq
\left(\int_{0}^{1}|\psi(x)|^{p}dx_{2}\right)^{\frac{1}{p}}\left(\int_{0}^{1}\frac{dx_{1}}{|1+\lambda x_{2}^{3}+\lambda x_{2}x_{1}^{2}|^{q}}\right)^{\frac{1}{q}}.$$

Then using again the  H\"older   inequality for  this integral we establish

$$|I_{\alpha,\beta}|\leq
\left(\int_{0}^{1}\int_{0}^{1}|\psi(x)|^{p}dx_{2}dx_{1}\right)^{\frac{1}{p}}\left(\int_{0}^{1}\int_{0}^{1}\frac{dx_{1}}{|1+\lambda x_{2}^{3}+\lambda x_{2}x_{1}^{2}|^{q}}dx_{2}\right)^{\frac{1}{q}}.$$

Changing the variables as $x_{1}=\left(\frac{1+\lambda x_{2}^{3}}{\lambda x_{2}}\right)^{\frac{1}{2}}t$ we get

$$|I_{\alpha,\beta}|\leq\|\psi\|_{L^{p}}\left(\int_{0}^{1}\int_{0}^{1}\frac{dx_{1}dx_{2}}{|1+\lambda x_{2}^{3}+\lambda x_{2}x_{1}^{2}|^{q}}\right)^{\frac{1}{q}}$$$$=\|\psi\|_{L^{p}}\left(\int_{0}^{1}\frac{(1+\lambda x_{2}^{3})^{\frac{1}{2}-q}}{(\lambda x_{2})^{\frac{1}{2}}}dx_{2}\int_{0}^{A}\frac{dt}{(1+t^{2})^{q}}\right)^{\frac{1}{q}},$$
where $A=\left(\frac{\lambda x_{2}}{1+\lambda x_{2}^{3}}\right)^{\frac{1}{2}}$ and $\int_{0}^{A}\frac{dt}{(1+t^{2})^{q}}<C$ as $A\rightarrow\infty.$
Thus,
$$|I_{\alpha,\beta}|\leq C\|\psi\|_{L^{p}}\left(\int_{0}^{1}\frac{(1+\lambda x_{2}^{3})^{\frac{1}{2}-q}}{(\lambda x_{2})^{\frac{1}{2}}}dx_{2}\right)^{\frac{1}{q}}.$$
Replacing $x_{2}$ by $\lambda^{-\frac{1}{3}}\tau$ and using $\frac{1}{q}=1-\frac{1}{p}$ we get

$$|I_{\alpha,\beta}|\leq\frac{C\|\psi\|_{L^{p}}}{\lambda^{\frac{2}{3}-\frac{1}{3p}}}
\left(\int_{0}^{\lambda^{\frac{1}{3}}}\frac{d\tau}{\tau^{\frac{1}{2}}(\tau^{3}+1)^{q-\frac{1}{2}}}\right)^{\frac{1}{q}}\leq
\frac{C\|\psi\|_{L^{p}}}{\lambda^{\frac{2}{3}-\frac{1}{3p}}}
\left(\int_{0}^{\infty}\frac{d\tau}{\tau^{\frac{1}{2}}(\tau^{3}+1)^{q-\frac{1}{2}}}\right)^{\frac{1}{q}}.$$
Since the last integral is covergent,
$$|I_{\alpha,\beta}|\leq\frac{C\|\psi\|_{L^{p}}}{\lambda^{\frac{2}{3}-\frac{1}{3p}}}.$$

\textbf{Case III.} Assume that the phase function has critical point of type $A_{2}$ so that
$f(x)=x_{1}^{3}.$
We estimate the integral \eqref{int1} with phase function $f(x)=x_{1}^{3}$

$$|I_{\alpha,\beta}|\leq\int_{0}^{1}\int_{0}^{1}|E_{\alpha,\beta}(i\lambda x_{1}^{3})||\psi(x)|dx_{1}dx_{2}.$$
First we use inequality \eqref{int10} for the last inner integral we obtain

$$|J_{in3}|:=\int_{0}^{1}\frac{|\psi(x)|dx_{1}}{1+\lambda x_{1}^{3}}.$$
Then we use H\"older inequality for the last integral  $I_{\alpha,\beta}$ twise and we get:
$$|I_{\alpha,\beta}|\leq\left(\int_{0}^{1}\int_{0}^{1}|\psi(x)|^{p}dx_{1}dx_{2}\right)^{\frac{1}{p}}\left(\int_{0}^{1}\int_{0}^{1}\frac{dx_{1}}{|1+\lambda x_{1}^{3}|^{q}}dx_{2}\right)^{\frac{1}{q}}.$$
Replacing $\lambda^{-\frac{1}{3}}x_{1}$ by $t$ in the above inequality, we obtain
$$|I_{\alpha,\beta}|\leq\frac{C\|\psi\|_{L^{p}}}{\lambda^{\frac{1}{3q}}}\left(\int_{0}^{\lambda^{\frac{1}{3}}}\frac{dt}{|1+ t^{3}|^{q}}\right)^{\frac{1}{q}}\leq\frac{C\|\psi\|_{L^{p}}}{\lambda^{\frac{1}{3q}}}\left(\int_{0}^{\infty}\frac{dt}{|1+ t^{3}|^{q}}\right)^{\frac{1}{q}}.$$
Since $\frac{1}{p}+\frac{1}{q}=1$ and the last integral converges,
$$|I_{\alpha,\beta}|\leq\frac{C\|\psi\|_{L^{p}}}{\lambda^{\frac{1}{3}-\frac{1}{3p}}}.$$
The proof is complete.

\textbf{Remark.}
If $\alpha=\beta=1$ in the integral \eqref{int1} integral called oscillatory integral and theorem holds for it.

\textbf{Declaration of competing interest}

This work does not have any conflicts of interest.

\textbf{Acknowledgements.} This paper was supported by "El-yurt umidi" Foundation of Uzbekistan and partially supported
in parts by the FWO Odysseus 1 grant G.0H94.18N: Analysis and Partial Differential Equations and by the Methusalem programme of the Ghent University Special Research Fund (BOF) (Grant number 01M01021). The author thanks to Prof. M.Ruzhansky for
proposing the problem and constant attention to the work and also thanks
the referees for numerous suggestions which greatly helped to improve the
exposition.

\end{document}